\documentclass[12pt]{article}
\begin{document}
\title{
 A generalized Numerov method\\
 for linear second-order
 differential equations\\
 involving a first derivative term}
\author{%
V. I. Tselyaev \\
{\it \small Nuclear Physics Department,
 V. A. Fock Institute of Physics},\\
{\it \small St. Petersburg State University, 198504,
 St. Petersburg, Russia}}
\date{}
\maketitle
\begin{abstract}
The Numerov method for linear second-order differential
equations is generalized to include equations containing a first
derivative term. The method presented has the same degree
of accuracy as the ordinary Numerov sixth-order method.
A general scheme of the application to the numerical solution
of the Hartree-Fock equations is considered.
\end{abstract}

\vspace{2em}
\newpage

\section{Introduction}

The linear second-order differential equations of the type
 \begin{equation}
y''(x) + g(x)\,y'(x) + f(x)\,y(x) = 0
 \label{1}
 \end{equation}
occur in various fields of physics. Here we shall mean the problem
of the numerical self-consistent solution of the Hartree-Fock (HF)
equations of motion deduced from the Skyrme energy functional
\cite{VB} describing ground-state properties of atomic nuclei.
For spherical nuclei these equations can be reduced to the
form~(\ref{1}), the first derivative term arising due to the
radial dependence of the nucleon effective mass. Usually the task
is solved by the Runge-Kutta method. But this method
is not the best one for the HF self-consistency procedure
because it requires the interpolation of functions
$g(x)$ and $f(x)$ in Eq.~(\ref{1}) between the grid points
where function $y(x)$ is not calculated. Another well-known
method (see, for example, Refs.~\cite{H,B}),
which was proposed by B.~V.~Numerov in 1923,
is free of the pointed difficulty but the first derivative term
in Eq.~(\ref{1}) precludes from its immediate application.
Several modifications of the Numerov method (NM) were developed
\cite{LW,DMB} in order to include equations of the type~(\ref{1})
and more general nonlinear equations. Here another generalization
of the NM is presented which is most suitable for the
HF calculations \cite{T} and yields the same degree of accuracy
as the initial Numerov method.

\section{Generalized linear Numerov method}

     Let us introduce notations: $y_0 = y(x_0)$,
$y_{\pm} = y(x_0 \pm h)$ and analogously for $g(x)$ and $f(x)$,
where $x_0$ is the fixed grid point, $h$ is the step length.
Developing quantities $y_{\pm}$, $y'_{\pm}$, $y''_{\pm}$
in powers of $h$, we obtain
 \begin{eqnarray}
y_+ + y_- & = &
2 y_0 + h^2 y''_0 + \frac {h^4}{12} y^{(4)}_0
+ \frac {h^6}{360} y^{(6)}_0 + O (h^8)\,,
\label{2}\\
y_+ - y_- & = &
2 h y'_0 + \frac {h^3}{3} y'''_0
+ \frac {h^5}{60} y^{(5)}_0 + O (h^7)\,,
\label{s1}\\
y'_+ + y'_- & = &
2 y'_0 + h^2 y'''_0 + \frac {h^4}{12} y^{(5)}_0
+ O (h^6)\,,
\label{s2}\\
y'_+ - y'_- & = &
2 h y''_0 + \frac {h^3}{3} y^{(4)}_0
+ \frac {h^5}{60} y^{(6)}_0 + O (h^7)\,,
\label{s3}\\
y''_+ + y''_- & = &
2 y''_0 + h^2 y^{(4)}_0 + \frac {h^4}{12} y^{(6)}_0
+ O (h^6)\,,
\label{s4}\\
y''_+ - y''_- & = &
2 h y'''_0 + \frac {h^3}{3} y^{(5)}_0 + O (h^5)\,.
\label{s5}
 \end{eqnarray}
In addition, Eq.~(\ref{1}) yields:
 \begin{eqnarray}
y''_0 & = & - g_0 y'_0 - f_0 y_0\,,
\label{s6}\\
y''_+ & = & - g_+ y'_+ - f_+ y_+\,,
\label{s7}\\
y''_- & = & - g_- y'_- - f_- y_-\,.
\label{s8}
 \end{eqnarray}
Let us consider Eqs.~(\ref{s1})--(\ref{s8}) as a system
of eight linear equations for eight unknown quantities:
$y'_0$, $y''_0$, $y'''_0$, $y^{(4)}_0$, $y'_{\pm}$, $y''_{\pm}$.
Solving these equations and substituting the found quantities
$y''_0$ and $y^{(4)}_0$ in Eq.~(\ref{2}) we get
 \begin{equation}
T_0 y_0 = T_+ y_+ + T_- y_- + R^{(6)}\,,
 \label{3}
 \end{equation}
where
 \begin{eqnarray}
T_0 & = & 2\, a - \frac {5 h^2}{6} b_0 f_0\,,
\label{y1}\\
T_{\pm} & = & a \pm \frac {h}{24} (10\, c\, g_0 + g_+ + g_-)
+ \frac {h^2}{12} b_{\pm} f_{\pm}\,,
\label{y2}\\
a & = & \left( 1 + \frac {h}{3} g_+ \right)
\left( 1 - \frac {h}{3} g_- \right)
+ \frac {h^2}{18} g_0 (g_+ + g_-)\,,
\label{y3}\\
b_0 & = & \left( 1 + \frac {4 h}{15} g_+ \right)
\left( 1 - \frac {4 h}{15} g_- \right)
+ \left( \frac {h}{15} \right)^2 g_+ g_-\,,
\label{y4}\\
b_{\pm} & = &
\left( 1 \pm \frac {5 h}{6} g_0 \right)
\left( 1 \mp \frac {h}{3} g_{\mp} \right)
+ \left( \frac {h}{3} \right)^2 g_0 g_{\mp}\,,
\label{y5}\\
c & = & \left( 1 + \frac {7 h}{20} g_+ \right)
\left( 1 - \frac {7 h}{20} g_- \right)
+ \left( \frac {3 h}{20} \right)^2 g_+ g_-\,,
\label{y6}\\
R^{(6)} & = & \frac {h^6}{240}
(y_0^{(6)} + 3 g_0 y_0^{(5)}) + O (h^8)\,.
\label{y7}
 \end{eqnarray}

     In more detail this result can be obtained by the
following way. Making use of Eqs.~(\ref{s6})--(\ref{s8}),
we get from Eqs.~(\ref{s5}) and (\ref{s4})
 \begin{eqnarray}
y'''_0 & = & \frac{g_- y'_- + f_- y_- - g_+ y'_+ - f_+ y_+}
{2h} - \frac{h^2}{6} y^{(5)}_0 + O (h^4)\,,
\label{1c}\\
y^{(4)}_0 & = &
\frac{2(g_0 y'_0 + f_0 y_0) - g_+ y'_+ - f_+ y_+
- g_- y'_- - f_- y_-}{h^2}
\nonumber\\
&&- \frac{h^2}{12} y^{(6)}_0 + O (h^4)\,.
\label{2c}
 \end{eqnarray}
Substituting these equalities into Eqs.~(\ref{s1})--(\ref{s3})
we obtain
 \begin{eqnarray}
y_+ - y_- & = & 2hy'_0 + \frac{h^2}{6} 
( g_- y'_- + f_- y_- - g_+ y'_+ - f_+ y_+ )
\nonumber\\
&&- \frac{7h^5}{180} y^{(5)}_0 + O (h^7)\,,
\label{3c}\\
y'_+ + y'_- & = & 2y'_0 + \frac{h}{2}
( g_- y'_- + f_- y_- - g_+ y'_+ - f_+ y_+ )
\nonumber\\
&&- \frac{h^4}{12} y^{(5)}_0 + O (h^6)\,,
\label{4c}\\
y'_+ - y'_- & = & - \frac{h}{3} ( 4 g_0 y'_0 + 4 f_0 y_0
+ g_+ y'_+ + f_+ y_+ + g_- y'_- + f_- y_-  )
\nonumber\\
&&- \frac{h^5}{90} y^{(6)}_0 + O (h^7)\,.
\label{5c}
 \end{eqnarray}
It is useful to rewrite Eq.~(\ref{3c}) in the form:
 \begin{eqnarray}
y'_0 & = & \frac{y_+ - y_-}{2h}
+ \frac{h}{12} ( g_+ y'_+ + f_+ y_+ - g_- y'_- - f_- y_-  )
\nonumber\\
&&+ \frac{7h^4}{360} y^{(5)}_0 + O (h^6)\,.
\label{6c}
 \end{eqnarray}
Substitution for this formula into  Eqs.~(\ref{4c}) and
(\ref{5c}) leads to the following system of two equations
for quantities $y'_{\pm}$
 \begin{eqnarray}
\alpha_+ y'_+ + \alpha_- y'_- & = & u\,,
\label{7c}\\
\beta_+  y'_+ - \beta_-  y'_- & = & v\,,
\label{8c}
 \end{eqnarray}
where
 \begin{eqnarray}
\alpha_{\pm} & = & 1 \pm \frac{h}{3} g_{\pm}\,, \qquad
\beta_{\pm} = \alpha_{\pm} + \frac{h^2}{9} g_0 g_{\pm}\,,
\label{9c}\\
u & = & \frac{1}{h} (y_+ - y_-)
+ \frac{h}{3} (f_- y_- - f_+ y_+)
- \frac{2h^4}{45} y^{(5)}_0 + O (h^6)\,,
\label{10c}\\
v & = & \frac{2}{3} g_0 (y_- - y_+)
- \frac{h}{3} (4 f_0 y_0 + f_+ y_+ + f_- y_-)
\nonumber\\
&&+ \frac{h^2}{9} g_0 (f_- y_- - f_+ y_+)
- \frac{7h^5}{270} g_0 y^{(5)}_0 
- \frac{h^5}{90} y^{(6)}_0 + O (h^7)\,.
\label{11c}
 \end{eqnarray}
The solution of the system (\ref{7c}), (\ref{8c}) is as follows:
 \begin{equation}
y'_{\pm} = \frac{\beta_{\mp} u \pm \alpha_{\mp} v}{2a}\,,
\label{12c}
 \end{equation}
where $a = \frac{1}{2} (\alpha_+ \beta_- + \alpha_- \beta_+)$,
that coincides with the definition (\ref{y3}).

     Substituting the formulas (\ref{12c}) into right-hand side
of Eq.~(\ref{6c}) we obtain after some algebra with taking into
account Eqs.~(\ref{9c})--(\ref{11c})
 \begin{equation}
y'_0 = \frac {S_+ y_+ - S_- y_- - S_0 y_0}{2 a h} + R^{(4)}\,,
 \label{4}
 \end{equation}
where
 \begin{eqnarray}
S_0 & = & \frac {h^3}{9} (g_+ + g_-) f_0\,,
\label{dy1}\\
S_{\pm} & = & \left( 1 + \frac {5 h}{12} g_+ \right)
\left( 1 - \frac {5 h}{12} g_- \right)
+ \left( \frac {h}{12} \right)^2 g_+ g_- + \frac {h^2}{6}
\left( 1 \mp \frac {h}{3} g_{\mp} \right) f_{\pm}\,,
\label{dy2}\\
R^{(4)} & = & \frac {7 h^4}{360} y_0^{(5)} + O (h^6)\,.
\label{dy3}
 \end{eqnarray}

     Finally, substitution for the found solutions
$y'_{\pm}$, $y'_0$ (Eqs.~(\ref{12c}), (\ref{4})) into
Eqs.~(\ref{s6}) and (\ref{2c}) yields the explicit formulas
for the quantities $y''_0$ and $y^{(4)}_0$ in terms of
$y_0$, $y_{\pm}$. After substitution for these formulas into
Eq.~(\ref{2}) and a series of lengthy but straightforward
algebraic transformations we arrive at the result
(\ref{3})--(\ref{y7}).
Omitting the term $R^{(6)}$ in Eq.~(\ref{3}) we obtain
the recurrence three-point formula of the generalized
Numerov method for linear second-order differential
equations or, for brevity, of the generalized linear NM (GLNM).
Clearly this method reduces to the ordinary NM if
$g(x) = 0$ in Eq.~(\ref{1}). The local truncation error
of the GLNM, which is contained in the term $R^{(6)}$, is one
of the same order $h^6$ as the error of the ordinary NM.

     The formula (\ref{4}) enables one to calculate the first
derivative if the function $y(x)$ is known at the grid points.
The local truncation error of order $h^4$ is determined by the
term $R^{(4)}$. The more precise formula follows immediately
from Eq.~(\ref{5c})
 \begin{equation}
y'_{\pm} = \frac {\left( 3 \mp h g_{\mp} \right)
y'_{\mp} \mp h \left( 4 g_0 y'_0
+ 4 f_0 y_0 + f_+ y_+ + f_- y_- \right) \mp 3 R^{(5)}}
{3 \pm h g_{\pm}}\,,
 \label{5}
 \end{equation}
where
 \begin{equation}
R^{(5)} = \frac {h^5}{90} y_0^{(6)}
+ O (h^7)\,.
 \label{6}
 \end{equation}
The presence of the derivatives in the right-hand side
is a shortcoming of this formula, nevertheless
Eq.~(\ref{5}) is practical for the evaluation of $y'(x)$
at the endpoints of the grid.

\section{Application to the Hartree-Fock calculations}

     Consider a general scheme within which the method proposed
can be applied to the numerical solution of the HF equations.
The HF approximation is a basis of numerous microscopic physical
theories describing the quantum many-body systems. In general
formulation, the HF method leads to a system of nonlinear
integrodifferential equations. We shall consider a special case
of the HF equations of motion deduced from the energy functional
constructed on the base of
the zero-range Skyrme forces \cite{VB} which are
widely used for the description of atomic nuclei properties
(see, for example, Ref.~\cite{LGQ}). The variational principle
applied to the Skyrme energy functional leads in the case of
spherical nuclei to the following system of equations
(in proper units):
 \begin{equation}
z''_{\lambda}(x) - \frac{m'_q (x)}{m_q (x)} z'_{\lambda}(x)
+ 2m_q (x) [e_{\lambda} - V_{\lambda}(x)] z_{\lambda}(x) = 0\,,
\label{1hf}
 \end{equation}
where $q$ denotes a sort of nucleon (proton or neutron),
the index $\lambda$ stands for the set of orbital quantum
numbers (including $q$), $x$ denotes the radial coordinate,
$z_{\lambda}(x)$ is the radial wave function, $e_{\lambda}$
is the eigenvalue playing the role of single-particle energy,
$V_{\lambda}(x)$ is the state-dependent mean-field potential,
$m_q (x)$ is the nucleon effective mass.

     Comparing Eqs.~(\ref{1}) and (\ref{1hf}), we see that
they have the same form and would be identical if we put
 \begin{equation}
y(x) = z_{\lambda}(x)\,, \quad
g(x) = - \frac{m'_q (x)}{m_q(x)}\,, \quad
f(x) = 2m_q (x) [e_{\lambda} - V_{\lambda}(x)]\,.
\label{2hf}
 \end{equation}
The essential difference consists in the following:
actually formula (\ref{1hf}) stands for the system of $N$
coupled nonlinear integrodifferential equations because
the quantities $m_q (x)$ and $V_{\lambda}(x)$ are functionals
of the densities which depend in turn
on the set $\{ z_{\lambda}, z'_{\lambda} \}$ of all the wave
functions and their derivatives with
$\lambda \in \{ \lambda_1, \ldots , \lambda_N \}$
(see \cite{VB}).
In practice, the system of equations (\ref{1hf}) is solved
by making use of some iteration procedure.
The convergence is achieved by averaging of the densities
calculated on two successive iterations.
The description and the analysis of the procedure in more
detail are outside the scope of the present paper (see, e. g.,
Ref.~\cite{SV}, and references therein where some relevant
methods are discussed). Here it is important only that on each
fixed HF iteration one has to solve the set of $N$
{\it uncoupled linear differential} equations which have
the same form (\ref{1hf}) but with already {\it known}
functions $m_q (x)$ and $V_{\lambda}(x)$ determined by the
results of previous iterations.

     The numerical integration of Eq.~(\ref{1hf}) in this case
can be performed by means of the GLNM described above.
For the sake of simplicity it is convenient to come back to the
notations of the preceding section taking into account
Eqs.~(\ref{2hf}). Setting
 \begin{eqnarray}
Y^{out}_0 & = & y_0 / y_+\,, \qquad
Y^{out}_- = y_- / y_0\,,
\label{3hf}\\
Y^{in}_0 & = & y_0 / y_-\,, \qquad
Y^{in}_+ = y_+ / y_0\,,
\label{4hf}
 \end{eqnarray}
we obtain the following recurrence relations from Eq.~(\ref{3})
(omitting the term $R^{(6)}$):
 \begin{eqnarray}
Y^{out}_0 & = & T_+ / (T_0 - T_- Y^{out}_-)\,,
\label{5hf}\\
Y^{in}_0  & = & T_- / (T_0 - T_+ Y^{in}_+)\,.
\label{6hf}
 \end{eqnarray}
Note that these transformations of Eq.~(\ref{3}) are similar
but not identical to ones of the renormalized Numerov method
developed in Ref.~\cite{J}.

In the proposed scheme Eqs.~(\ref{5hf}) and (\ref{6hf}) are
used for the outward and the inward integrations, respectively.
The outward integration starts from a point near $x=0$.
The initial value of $Y^{out}_-$ in Eq.~(\ref{5hf}) is
calculated using analytic expansions of the regular solutions
$z_{\lambda}(x)$ about $x=0$. To improve the accuracy of
calculations it is practical to decrease the step length of the
grid near this point.
The inward integration starts from an outside endpoint
of the grid where the irregular solutions of
Eq.~(\ref{1hf}) are known analytically. The initial value of
$Y^{in}_+$ in Eq.~(\ref{6hf}) is calculated as the ratio of the
Whittaker functions for protons and of the spherical Hankel ones
for neutrons. The eigenvalue $e_{\lambda}$ is found from the
condition $Y^{out}_0 = 1 / Y^{in}_+$ at some matching point
$x_0 = x_m$. The reasonable choice for this point is
(see, e. g., Ref.~\cite{VB})
the approximate position of the last extremum
of the wave function $z_{\lambda}(x)$. In the end of this
procedure, which is performed for each $z_{\lambda}(x)$
separately, the function $z_{\lambda}(x)$ is calculated at the
grid points using the ratios $Y^{out}_0$, $Y^{in}_0$, and the
normalization condition:
 \begin{equation}
\int^{\infty}_0 dx z^2_{\lambda}(x) = 1\,.
\label{7hf}
 \end{equation}

     Finally, the derivatives $z'_{\lambda}(x)$ are calculated
employing Eq.~(\ref{4}) at the inside points of the grid and
Eq.~(\ref{5}) at the endpoints. After this, new approximations
to the functions $m_q (x)$, $V_{\lambda}(x)$, which are used
in the next HF iteration, are calculated.

     The algorithm, that was briefly outlined above, has been
realized in the computer code \cite{T} intended for the
Skyrme-Hartree-Fock calculations. The convergence and stability
of the procedure described were tested in the calculations of
ground-state properties of all doubly magic atomic nuclei using
most of the present Skyrme-force parametrizations. It was
obtained that the algorithm based on the GLNM reproduce the
known reference results within their accuracy.

\newpage
\end{document}